\documentclass[a4paper,12pt]{article}
\usepackage[T1]{fontenc}
\usepackage[latin1]{inputenc}
\usepackage{amsthm,amsmath}
\usepackage{amssymb}
\usepackage[all]{xy}
\usepackage{enumerate}
\usepackage{a4wide}
\usepackage{hyperref}
\usepackage{latexsym}
\usepackage{enumerate}
\usepackage[mathscr]{eucal}

\newtheorem{theorem}{Theorem}[section]
\newtheorem{proposition}[theorem]{Proposition}
\newtheorem{corollary}[theorem]{Corollary}

\newtheorem*{theorem*}{Theorem}
\newtheorem*{proposition*}{Proposition}
\newtheorem*{corollary*}{Corollary}
\newtheorem*{lemma*}{Lemma}
\theoremstyle{definition}
\newtheorem{definition}[theorem]{Definition}

\newtheorem{example}[theorem]{Example}
\newtheorem{remark}[theorem]{Remark}
\newtheorem*{remark*}{Remark}

\newtheorem*{definition*}{Definition}
\numberwithin{equation}{section}

\newcommand{\cat}[1]{\mathcal{#1}}
\newcommand{\coring}[1]{\mathfrak{#1}}
\newcommand{\tensor}[1]{\otimes_{#1}}

\newcommand{\rcomod}[1]{ \mathsf{Comod}_{#1}}
\newcommand{\rmod}[1]{\mathsf{Mod}_{#1}}
\newcommand{\lmod}[1]{{}_{#1}\mathsf{Mod}}

\newcommand{\cotensor}[1]{\square_{#1}}

\newcommand{\homcom}[3]{\mathrm{Hom}_{#1}(#2,#3)}
\renewcommand{\hom}[3]{\mathrm{Hom}_{#1}\left( #2 \, , \, #3\right)}

\newcommand{\rend}[2]{\mathrm{End}({#2}_{#1})}
\newcommand{\lend}[2]{\mathrm{End}({}_{#1}#2)}

\newcommand{\rcomatrix}[2]{#2^* \tensor{#1} #2}

\newcommand{\entmod}[2]{\mathcal{M}(\Psi)_{#1}^{#2}}

\newcommand{\fk}[1]{\mathfrak{#1}}
\newcommand{\Sf}[1]{\mathsf{#1}}

\newcommand{\esc}[2]{\langle #1,#2 \rangle}

\newcommand{\Rat}{\mathrm{Rat}}

\newcommand{\can}[1]{\mathsf{can}_{#1}}

\newcommand{\lr}[1]{\left(\underset{}{} #1 \right)}

\newcommand{\TT}{\mathscr{T}}
\newcommand{\T}{\mathcal{T}}

\begin{document}
\title{Corings with decomposition and semiperfect corings
\footnote{This research is supported by the grants MTM2004-01406
and MTM2007-61673 from the Ministerio de Edu\-ca\-ci\'{o}n y Ciencia
of Spain, and P06-FQM-01889 from the Consejer\'{\i}a de Innovaci\'{o}n,
Ciencia y Empresa of Andaluc\'{\i}a, with funds from FEDER (Uni\'{o}n
Europea)} }
\author{L. El Kaoutit \\
\normalsize Departamento de \'{A}lgebra \\ \normalsize Facultad de Educaci\'{o}n y Humanidades \\
\normalsize  Universidad de Granada \\ \normalsize El Greco N$^o$
10, E-51002 Ceuta, Espa\~{n}a \\ \normalsize
e-mail:\textsf{ kaoutit@ugr.es} \and J. G\'omez-Torrecillas \\
\normalsize Departamento de \'{A}lgebra \\ \normalsize Facultad de Ciencias \\
\normalsize Universidad
de Granada\\ \normalsize E18071 Granada, Espa\~{n}a \\
\normalsize e-mail: \textsf{gomezj@ugr.es} }

\date{ }
\maketitle

\begin{center}
\textbf{\textit{Opgedragen aan Freddy Van Oystaeyen ter
gelegenheid van zijn zestigste verjaardag}}
\end{center}

\section*{Introduction}

When studying generalized module categories like Doi-Koppinen Hopf
modules or, more generally, entwined modules (see
\cite{Caenepeel/Militaru/Zhu:2002} or
\cite{Brzezinski/Wisbauer:2003} for detailed accounts of their
fundamental properties), it is reasonable, after the observation
that they are categories of comodules \cite{Brzezinski:2002}, to
formulate our questions in terms of comodules over corings.
However, it happens often that the answer is given in terms of
abstract categorical concepts, which do not say much more on the
concrete comodule algebra or entwining structure we are dealing
with. In a more favorable situation, we could characterize
categorical properties of the category of relative modules in
terms of structural properties of the associated coring. But,
sometimes, the structure of the corings is described in terms
which have no direct relationship with the entwined algebra and
coalgebra. One example of the sketched situation is the following:
assume we want to study under which conditions the category
$\entmod{A}{C}$ of entwined right modules over an entwining
structure $(A,C)_{\Psi}$ for an algebra $A$ and a coalgebra $C$
has a generating set of small projectives. Of course, when
$\entmod{A}{C}$ is a Grothendieck category, a first answer is
given by Freyd's theorem: it has to be equivalent to a category of
modules over a ring with enough idempotents. Looking at
$\entmod{A}{C}$ as the category of right comodules over a suitable
coring \cite{Brzezinski:2002}, we have a characterization in terms
of infinite comatrix corings in \cite[Theorem
2.7]{ElKaoutit/Gomez:2004b}. However, none of these results give
answers directly expressable in terms of $A$ or $C$. In this
paper, we describe a class of $A$--corings for which there exists
a generating set of small projective objects for their category of
right comodules. As a consequence, we will obtain that if $C$ is a
right semiperfect coalgebra over a field, then $\entmod{A}{C}$ has
a generating set of small projective objects for every entwining
structure $(A,C)_{\Psi}$ (Corollary \ref{entgenerators}). Of
course, we think that the interest of our results go beyond the
theory of entwined modules, as the categories of comodules over
corings deserve to be investigated by their own right.

The idea comes from the theory of coalgebras. It is well known
\cite{Lin:1977} that if $C$ is a coalgebra over a field, then $C$
is right semiperfect (or, equivalently, its category of right
comodules $\rcomod{C}$ has a generating set consisting of
finite-dimensional projectives) if and only if $C$ decomposes as a
direct sum of finite-dimensional left subcomodules. We give a
generalization of this characterization to the case of corings
over an arbitrary ring $A$ (Theorem \ref{TeoremaII}). It is then
natural to try to understand right semiperfect corings. This is
done in Section \ref{sperfect}. Our point of view here, in
contrast with \cite{ElKaoutit/Gomez:2005} and
\cite{Caenepeel/Iovanov:2005}, deliberatively avoids the
assumption of conditions on the ground ring $A$. Our approach
rests upon the study of general semiperfect categories due to
Harada \cite{Harada:1973a}, and on the study of corings having a
generating set of small projective comodules developed in
\cite{ElKaoutit/Gomez:2004b}. We also discuss the notion of a
(right) local coring and the exactness of the rational functor.

\medskip

\textbf{Notations and basic notions.} We work over a commutative
ground base ring with $1$ denoted by $K$. The letters $A$, $B$ are
reserved to denote associative $K$--algebras with unit, which will
referred to as rings. A module over a ring with unit means an
unital module, and all bimodules are assumed to be central
$K$--bimodules. The category of all right $A$--modules is denoted
by $\rmod{A}$. A linear morphism acts on the left, so some
conventions should be established. That is, if ${}_AN$ is a left
$A$--module then its endomorphism ring $\lend{A}{N}$ is considered
as ring with the opposite multiplication of the usual composition
law. In this way $N$ is an $(A,\lend{A}{N})$--bimodule. While, if
$N_A$ is a right $A$--module, then its endomorphism ring
$\rend{A}{N}$ has a multiplication the usual composition, and $N$
becomes obviously an $(\rend{A}{N},A)$--bimodule. For any
$(B,A)$--bimodule $M$, we consider in a canonical way its right
and left dual modules $M^*=\hom{}{M_A}{A_A}$,
${}^*M=\hom{}{{}_BM}{{}_BB}$ as $(A,B)$--bimodules.

We will also consider some rings without unit. When this is the
case will be clear by the context.

For any category $\cat{G}$, the notation $X \in \cat{G}$ means
that $X$ is an object of $\cat{G}$, and the identity morphism of
any object will be represented by the object itself.

Recall from \cite{Sweedler:1975} that an $A$--coring is a
three-tuple
$(\coring{C},\Delta_{\coring{C}},\varepsilon_{\coring{C}})$
consisting of an $A$--bimodule $\coring{C}$ and two homomorphisms
of $A$--bimodules
$$\xymatrix@C=50pt{\coring{C} \ar@{->}^-{\Delta_{\coring{C}}}[r] &
\coring{C}\tensor{A}\coring{C}},\quad \xymatrix@C=30pt{ \coring{C}
\ar@{->}^-{\varepsilon_{\coring{C}}}[r] & A}$$ such that
$(\Delta_{\coring{C}}\tensor{A}\coring{C}) \circ
\Delta_{\coring{C}} = (\coring{C}\tensor{A}\Delta_{\coring{C}})
\circ \Delta_{\coring{C}}$ and
$(\varepsilon_{\coring{C}}\tensor{A}\coring{C}) \circ
\Delta_{\coring{C}}=(\coring{C}\tensor{A}\varepsilon_{\coring{C}})
\circ \Delta_{\coring{C}}= \coring{C}$. A homomorphism of
$A$--corings is an $A$--bilinear map $\phi: \coring{C} \rightarrow
\coring{C}'$ which satisfies $\varepsilon_{\coring{C}'} \circ \phi
= \varepsilon_{\coring{C}}$ and $ \Delta_{\coring{C}'} \circ \phi
= (\phi \tensor{A} \phi) \circ \Delta_{\coring{C}}$.

A right $\coring{C}$--comodule is a pair $(M,\rho_{M})$ consisting
of a right $A$--module $M$ and a right $A$--linear map $\rho_{M}:
M \rightarrow M\tensor{A}\coring{C}$, called right
$\coring{C}$--coaction, such that
$(M\tensor{A}\Delta_{\coring{C}}) \circ \rho_M =
(\rho_M\tensor{A}\coring{C}) \circ \rho_M$ and
$(M\tensor{A}\varepsilon_{\coring{C}}) \circ \rho_M=M$. A morphism
of right $\coring{C}$--comodules is a right $A$--linear map $f: M
\rightarrow M'$ satisfying $ \rho_{M'} \circ f =
(f\tensor{A}\coring{C}) \circ \rho_M$. The $K$--module of all
homomorphisms of right $\coring{C}$--comodules from a comodule
$M_{\coring{C}}$ to a comodule $M'_{\coring{C}}$ is denoted by
$\homcom{\coring{C}}{M}{M'}$. Right $\coring{C}$--comodules and
their morphisms form a $K$--linear category $\rcomod{\coring{C}}$
which is a Grothendieck category provided ${}_A\coring{C}$ is a
flat module, see \cite[Section 1]{ElKaoutit/Gomez/Lobillo:2004c}.
Left $\coring{C}$--comodules and their morphisms are symmetrically
defined. If $P$ is a right $\coring{C}$--comodule such that $P_A$
is finitely generated and projective (\emph{profinite}, for
short), then its right dual $P^*$ admits, in a natural way, a
structure of left $\coring{C}$--comodule
\cite[19.19]{Brzezinski/Wisbauer:2003}. The same arguments are
pertinent for left $\coring{C}$--comodules. The natural
isomorphism $P \cong {}^*(P^*)$ of right $A$--modules is in fact
an isomorphism of $\coring{C}$--comodules. Furthermore, if $P$ and
$Q$ are two right $\coring{C}$--comodules profinite as right
$A$--modules, then the right dual functor induces a $K$--module
isomorphism $\homcom{\coring{C}}{P}{Q} \cong
\homcom{\coring{C}}{Q^*}{P^*}$.

\section{Comatrix corings and corings with decompositions}\label{decompsition}

Let $\cat{P}$ be a set of profinite modules over a ring with unit
$A$. Consider $\Sigma = \oplus_{P \in \cat{P}} P$ and, for each $P
\in \cat{P}$, let $\iota_P : P \rightarrow \Sigma$, $\pi_P :
\Sigma \rightarrow P$ be, respectively, the canonical inclusion
and the canonical projection. Consider the set $\{ u_P = \iota_P
\circ \pi_P : P \in \cat{P} \}$ of ortogonal idempotents of
$\rend{A}{\Sigma}$, and let $T$ be a unital subring of
$\rend{A}{\Sigma}$ that contains the idempotents $u_P$. Write
\begin{equation*}
R = \bigoplus_{P,Q} u_QTu_P,
\end{equation*}
a ring with enough idempotents. Its category of right unital
$R$--modules is
 denoted by $\rmod{R}$. Unital here means $MR = M$, for $M$ a right $R$--module.
 Let us recall one of the three
constructions given in \cite{ElKaoutit/Gomez:2004b} for the
(infinite) comatrix $A$--coring associated to $\cat{P}$ an $T$ (or
$R$). We have now the $(R,A)$--bimodule $\Sigma$ and the
$(A,R)$--bimodule
 $\Sigma^{\dag}= \oplus_{P \in \cat{A}}
P^*$. Both $\Sigma^{\dag}$ and $\Sigma$ are unital $R$--modules.
In the $A$--bimodule $\Sigma^{\dag} \tensor{R} \Sigma$ we have
that $\phi \tensor{R} x = 0$ whenever $\phi \in P^*, x \in Q$ for
$P \neq Q$. In this way, the formula
\[
\Delta (\phi \tensor{R} x) = \sum \phi \tensor{R} e_{P,i}
\tensor{A} e_{P,i}^* \tensor{R} x \qquad (\phi \in P^*, \; x \in
P),
\]
where $\{ e_{P,i},e_{P,i}^* \}$ is a dual basis for the profinite
right $A$--module $P$, determines a comultiplication
\[
\xymatrix{\Delta : \Sigma^{\dag} \tensor{R} \Sigma \ar[rr] & &
\Sigma^{\dag} \tensor{R} \Sigma \tensor{A} \Sigma^{\dag}
\tensor{R} \Sigma}
\]
which, according to \cite[Proposition 5.2]{ElKaoutit/Gomez:2004b},
endows $\Sigma^{\dag} \tensor{R} \Sigma$ with a structure of an
$A$--coring referred to as the \emph{infinite comatrix coring}
associated to the set $\cat{P}$ and the ring $R$. Its counit is
given by evaluation of the forms $\phi \in P^*$ at the elements $x
\in P$, when $P$ runs $\cat{P}$.

When $\cat{P}$ consists of right comodules over an $A$--coring
$\coring{C}$, we have that $\Sigma = \bigoplus_{P \in \cat{P}}P$
is a right $\coring{C}$--comodule. We consider then the infinite
comatrix coring by putting $T = \rend{\coring{C}}{\Sigma}$. We
have that

\begin{equation}\label{sufidem}
R = \bigoplus_{P,Q \in \cat{P}}u_Q\rend{\coring{C}}{\Sigma}u_P
\cong \bigoplus_{P,Q \in \cat{P}}\hom{\coring{C}}{P}{Q}
\end{equation}

It follows from \cite[Lemma 4.7]{ElKaoutit/Gomez:2004b} and
\cite[diagram (5.12)]{ElKaoutit/Gomez:2004b} that there is a
canonical homomorphism of $A$--corings
\begin{equation}\label{can}
\xymatrix@R=0pt{\can{}:
\Sigma^{\dag} \tensor{R}\Sigma \ar@{->}[r] & \coring{C} \\
\varphi \tensor{R} x \ar@{|->}[r] & (\varphi \tensor{A}
\coring{C})  \rho_{\Sigma}(x)},
\end{equation}
where $\varphi \in \Sigma^{\dag}$ acts on $y \in \Sigma$ by
evaluation in the obvious way.

\begin{definition}
The comodule $\Sigma$ is said to be $R-\coring{C}$--\emph{Galois}
if the canonical map $\mathsf{can}$ is bijective.
\end{definition}

Recall from \cite[Section 2]{ElKaoutit/Gomez/Lobillo:2004c} that a
three-tuple $\TT=(\coring{C},B,\esc{-}{-})$ consisting of an
$A$-coring $\coring{C}$, an $A$-ring $B$ (i.e., $B$ is an algebra
extension of $A$) and a balanced $A$--bilinear form $\esc{-}{-}:
\coring{C} \times B \to A$, is said to be \emph{a right rational
pairing over} $A$ provided
\begin{enumerate}[(1)]
 \item $\beta_A : B \to {}^*\coring{C}$ is a ring
     anti-homomorphism, where ${}^*\coring{C}$ is the left
     dual convolution ring of $\coring{C}$ defined in
     \cite[Proposition 3.2]{Sweedler:1975}, and

\item $\alpha_M$ is an injective map, for each right
    $A$-module $M$,
\end{enumerate}

where $\alpha_{-}$ and $\beta_{-}$ are the following  natural
transformations
$$\xymatrix@R=0pt{\beta_N : B \tensor{A} N \ar@{->}[r] &
\hom{}{_{A}\coring{C}}{{}_{A}N}, \\ b \tensor{A} n \ar@{->}[r]& \left[ c
\mapsto \esc{c}{b}n \right] }\quad \xymatrix@R=0pt{ \alpha_M : M
\tensor{A} \coring{C} \ar@{->}[r] & \hom{}{B_{A}}{M_{A}} \\
m \tensor{A} c \ar@{->}[r] & \left[ b \mapsto m\esc{c}{b} \right].
}$$

Given a right rational pairing $\TT=(\coring{C},B,\esc{-}{-})$
over $A$, we can define a functor called the \emph{right rational
functor} as follows.  An element $m$ of a right $B$--module $M$ is
called \emph{rational} if there exists a set of \emph{right
rational parameters} $\{(c_i,m_i)\} \subseteq \coring{C} \times M$
such that $ m b = \sum_i m_i\esc{c_i}{b}$, for all $b \in B$. The
set of all rational elements in $M$ is denoted by $\Rat^{\TT}(M)$.
As it was explained in \cite[Section
2]{ElKaoutit/Gomez/Lobillo:2004c}, the proofs detailed in
\cite[Section 2]{Gomez:1998} can be adapted in a straightforward
way in order to get that $\Rat^{\TT}(M)$ is a $B$--submodule of
$M$ and the assignment $M \mapsto \Rat^{\TT}(M)$ is a well defined
functor
\begin{equation*}
\Rat^{\TT} : \rmod{B} \rightarrow \rmod{B},
\end{equation*}
which is in fact a left exact preradical \cite[Ch.
VI]{Stenstrom:1975}. Therefore, the full subcategory
$\Rat^{\TT}(\rmod{B})$ of $\rmod{B}$ whose objects are those
$B$--modules $M$ such that $\Rat^{\TT}(M) = M$ is a closed
subcategory. Furthermore, $\Rat^{\TT}(\rmod{B})$ is a Grothendieck
category which is shown to be isomorphic to the category of right
comodules $\rcomod{\coring{C}}$ as \cite[Theorem
2.6']{ElKaoutit/Gomez/Lobillo:2004c} asserts (see also
\cite[Proposition 2.8]{Abuhlail:2003a}). We say that a set
$\cat{S}$ of objects of a Grothendieck category $\cat{A}$ is a
\emph{generating set} of $\cat{A}$ if the coproduct $\bigoplus_{X
\in \cat{S}} X$ is a generator of $\cat{A}$.

\medskip

We are ready to state and prove our main theorem.

\begin{theorem}\label{TeoremaII}
The following statements are equivalent for an $A$-coring
$\coring{C}$:
\begin{enumerate}[(i)]
\item\label{descomposicion} $\coring{C} = \bigoplus_{E \in
\cat{E}} E$, for a family of subcomodules $\cat{E}$ of
${}_{\coring{C}}\coring{C}$ such that ${}_AE$ is profinite for
every $E \in \cat{E}$; %
\item\label{categorico} ${}_A\coring{C}$ is projective and
    there exists a generating set $\cat{P}$ of small
    projective objects in $\rcomod{\coring{C}}$ such that the
    right comodule $\Sigma = \bigoplus_{P \in \cat{P}}P$,
    considered as a left module over $R = \bigoplus_{P,Q \in
    \cat{P}}\homcom{\coring{C}}{P}{Q}$, admits a decomposition
    as direct sum of finitely generated $R$-submodules;
\item\label{A-modulo} ${}_A\coring{C}$ is projective and there
    exists a generating set $\cat{P}$ of
    $\rcomod{\coring{C}}$, whose members are profinite as
    right $A$-modules, such that $\Sigma = \bigoplus_{P \in
    \cat{P}}P$ admits, as a left module over $R =
    \bigoplus_{P,Q \in \cat{P}}\homcom{\coring{C}}{P}{Q}$, a
    decomposition as direct sum of finitely generated
    $R$-submodules;
    \item\label{Galois} ${}_A\coring{C}$ is projective and
    there exists a set $\cat{P}$ of right
    $\coring{C}$-comodules profinite as right $A$-modules such
    that $\coring{C}$ is $R-\Sigma$-Galois for $\Sigma =
    \bigoplus_{P \in \cat{P}}P$ and $R = \bigoplus_{P,Q \in
    \cat{P}}\homcom{\coring{C}}{P}{Q}$, and $\Sigma$ admits,
    as a left $R$-module, a decomposition as direct sum of
    finitely generated $R$-submodules.
\end{enumerate}
\end{theorem}
\begin{proof}
 $(i) \Rightarrow (ii)$ Associated to the given decomposition of left comodules
 $\coring{C}=\oplus_{E \in\cat{E}}E$, there is a family of orthogonal idempotents
 $\{e_{E}:\,\, E \in \cat{E}\}$ in $\lend{\coring{C}}{\coring{C}}$,
 where $e_{E}= \iota_{E} \circ \pi_E$, for $\iota_E:E \to \coring{C}$ the
 canonical injection and $\pi_E: \coring{C} \to E$ the canonical projection for each $E \in \cat{E}$.
 The ring $\lend{\coring{C}}{\coring{C}}$ is endowed with the multiplication opposite to the composition.
 Since ${}_A\coring{C}$ is projective we have the canonical rational pairing
 $\TT=(\coring{C},\lend{\coring{C}}{\coring{C}},\esc{-}{-})$,
 where $\esc{c}{f}=\varepsilon(f(c))$, for $c \in \coring{C}$, $f \in \lend{\coring{C}}{\coring{C}}$.
 Thus each right $\coring{C}$-comodule admits a right $\lend{\coring{C}}{\coring{C}}$--action,
 and so is in particular for the right $\coring{C}$-comodules ${}^*E$, with $E \in \cat{E}$.

Consider the set of right $\coring{C}$-comodules
$\cat{P}=\{{}^*E:\, E \in \cat{E}\}$ and the right
$\coring{C}$-comodule $\Sigma = \oplus_{E \in \cat{E}} {}^*E$.
Each of the maps
\begin{equation}\label{Eq-Iso}
\xymatrix@R=0pt{ e_{E}\lend{\coring{C}}{\coring{C}}  \ar@{->}[r]
& \homcom{\coring{C}}{E}{\coring{C}} \cong {}^*E \\ e_{E}f \ar@{|->}[r]
& \varepsilon  \circ f \circ e_E \circ \iota_E }
\end{equation}
is an isomorphism of right
$\lend{\coring{C}}{\coring{C}}$-modules, which means that each
$e_{E}\lend{\coring{C}}{\coring{C}}$ is actually a rational right
$\lend{\coring{C}}{\coring{C}}$-module, and so a right
$\coring{C}$-comodule. In this way, we get an isomorphism of right
$\coring{C}$-comodules
\begin{eqnarray}\label{Eq-Sigma}
 \bigoplus_{E \in \cat{E}} e_{E}\lend{\coring{C}}{\coring{C}} \,\, \cong\,\,
 \bigoplus_{E \in \cat{E}} {}^*E =\Sigma.
\end{eqnarray}
Since $e_{E}\lend{\coring{C}}{\coring{C}}$ is a small object in
the category $\rcomod{\coring{C}}$ for every $E \in \cat{E}$, we
deduce that
\begin{eqnarray}\label{Eq-EF}
 e_{E}\lend{\coring{C}}{\coring{C}} \,\,=\,\, \bigoplus_{F \in\, \cat{E}} e_{E}\lend{\coring{C}}{\coring{C}}e_F.
\end{eqnarray}

Consider in the ring $\rend{\coring{C}}{\Sigma}$ (with
multiplication the usual composition), the set of idempotents
$u_{{}^*E}\,=\, \iota_{{}^*E} \circ \pi_{{}^*E}$, for
$\iota_{{}^*E}: {}^*E \to \Sigma$ the canonical injection and
$\pi_{{}^*E}: \Sigma \to {}^*E$ the canonical projection for each
$E \in \cat{E}$. We have the ring with enough orthogonal
idempotents
$$  R\,\,=\,\,\bigoplus_{E,\,F\, \in \,\cat{E}} u_{{}^*E}
\rend{\coring{C}}{\Sigma} u_{{}^*F}.$$ Clearly we already have,
for every pair of comodules $E$, $F \in \cat{E}$, $K$-linear
isomorphisms
\begin{eqnarray}\label{Eq-ECF}
 e_{E} \lend{\coring{C}}{\coring{C}} e_{F} \,\,\cong \,\, \homcom{\coring{C}}{E}{F},
\end{eqnarray}
\begin{eqnarray}\label{Eq-UEF}
 u_{{}^*E} \rend{\coring{C}}{\Sigma} u_{{}^*F} \,\,\cong \,\,
 \homcom{\coring{C}}{{}^*F}{{}^*E}.
\end{eqnarray}
Using equations \eqref{Eq-ECF} and \eqref{Eq-UEF} we define, taking into account the canonical $K$-linear isomorphism
$$\homcom{\coring{C}}{E}{F}\,\,\cong \,\, \homcom{\coring{C}}{{}^*F}{{}^*E},$$ an isomorphism
\begin{eqnarray}\label{Eq-EFFE}
 e_{E} \lend{\coring{C}}{\coring{C}} e_{F} \,\,\cong \,\, u_{{}^*E} \rend{\coring{C}}{\Sigma} u_{{}^*F}
\end{eqnarray}
for each pair $E,F \in \cat{E}$.

In view of equality \eqref{Eq-EF} and the family of isomorphisms \eqref{Eq-EFFE}, we deduce a $K$-linear isomorphism
\begin{eqnarray}\label{Eq-EFR}
\bigoplus_{E\,\in\,\cat{E}}e_E \lend{\coring{C}}{\coring{C}}
 \,\, \cong \,\, \bigoplus_{E\,\in\,\cat{E}} u_{{}^*E}
\rend{\coring{C}}{\Sigma} u_{{}^*F} =R
\end{eqnarray}
Now, if we compose the isomorphisms given in \eqref{Eq-Sigma} and \eqref{Eq-EFR} we get an isomorphism
\begin{eqnarray}\label{Eq-f}
 f:\Sigma \longrightarrow R.
\end{eqnarray}
In order to give an explicit expression for the isomorphism
\eqref{Eq-f}, we should take into account that the inverse map of
\eqref{Eq-Sigma} is defined as follows: to each element
$\theta_{E} \in {}^*E$ it corresponds
$(\coring{C}\tensor{A}\theta_E) \circ \Delta \circ e_E \in
e_E\lend{\coring{C}}{\coring{C}}$. From this, given
$\theta=\sum_{E\,\in\, \cat{E}} \theta_E \in \oplus_{E
\,\in\,\cat{E}}{}^*E =\Sigma$, we get a map
$$ f(\theta): \, \bigoplus_{F\,\in\,\cat{E}}{}^*F \longrightarrow \bigoplus_{F\,\in\,\cat{E}}{}^*F$$
defined by
\begin{eqnarray}\label{Eq-Expresion}
f(\theta) \lr{\sum_{F\,\in\,\cat{E}}\varphi_F} \,=\, \sum_{E,\,F\,
\in\,\cat{E}} \varphi_E \circ \pi_E \circ (\coring{C}
\tensor{A}\theta_E) \circ \Delta \circ \iota_E, \quad
\lr{\sum_{F\,\in\,\cat{E}}\varphi_F \,\in\,
\bigoplus_{F\,\in\,\cat{E}}{}^*F}
\end{eqnarray}
Using the expression \eqref{Eq-Expresion}  we can easily show that
$f:\Sigma \to R$ is in fact a left $R$-module isomorphism. In
particular, we deduce from the decomposition
$R=\oplus_{F\,\in\,\cat{E}} Ru_{{}^*F}$ that ${}_R\Sigma$ admits a
decomposition as direct sum of finitely generated $R$-submodules.

Consider now the infinite comatrix coring $\Sigma^{\dagger}
\tensor{R}\Sigma$, for $\Sigma^{\dagger}=\oplus_{
E\,\in\,\cat{E}}({}^*E)^*$. We have in fact an isomorphism
$\Sigma^{\dagger} \cong \coring{C}$ of left comodules. A routine
computation show that the following diagram
$$\xymatrix@R=30pt@C=60pt{ \Sigma^{\dagger} \tensor{R}\Sigma
 \ar@{->}^-{\Sf{can}}[r] \ar@{->}_-{\Sigma^{\dagger}\tensor{R}f}[d] & \coring{C}
\\ \Sigma^{\dagger} \tensor{R}R \ar@{->}^-{\cong}[r] & \Sigma^{\dagger} \ar@{->}_-{\cong}[u] } $$
is commutative, and so $\Sf{can}$ is an isomorphism. Since
${}_R\Sigma$  is a faithfully flat module because ${}_R\Sigma
\cong {}_RR$, we can apply \cite[Theorem 5.7($(iii)\Rightarrow
(i)$)]{ElKaoutit/Gomez:2004b} to deduce that $\{{}^*E:\, E \in
\cat{E}\}$ is a generating set of a small projectives for
$\rcomod{\coring{C}}$. \\
$(ii) \Rightarrow (iii)$ This clear, since each $P \in\cat{P}$ is,
as right $A$-module, finitely generated and projective (see
\cite[Theorem 5.7]{ElKaoutit/Gomez:2004b}).\\
 $(iii) \Rightarrow
(iv)$ This is deduced from \cite[Theorem 4.8]{ElKaoutit/Gomez:2004b}
and \cite[diagram (5.12)]{ElKaoutit/Gomez:2004b}.\\%
$(iv) \Rightarrow (i)$ We can consider $R$ as the (no unital)
subring of $\rend{\coring{C}}{\Sigma}$ given by $$R\,\,=\,\,
\bigoplus_{P,\,Q\, \in \, \cat{E}} u_P
\rend{\coring{C}}{\Sigma}u_Q,$$ where $u_P=\iota_P \circ \pi _P$
for $\pi_P: \Sigma \to P$ (resp. $\iota_P: P \to \Sigma$) is the
canonical projection (resp. injection). Consider the
decomposition ${}_R\Sigma=\oplus_{i \,\in\, I}\Sigma_i$ as direct
sum of finitely generated $R$-submodules ${}_R\Sigma_i$. We have
the following decomposition of the infinite comatrix coring as a
left comodule
\begin{eqnarray}\label{Eq-comatrix}
\Sigma^{\dagger} \tensor{R} \Sigma \,\, \cong \,\, \bigoplus_{i\, \in\, I}\Sigma^{\dagger}\tensor{R}\Sigma_i
\end{eqnarray}
For each $i \in I$ there exists a presentation of the left
$R$-module $\Sigma_i$
\begin{eqnarray}\label{Eq-represent}
F_i \longrightarrow \Sigma_i \longrightarrow 0,
\end{eqnarray}
where $F_i=\oplus_{P\,\in\, \cat{P}_i} Ru_p$, with a finite subset
$\cat{P}_i \subset \cat{P}$. Applying the functor
$\Sigma^{\dag}\tensor{R}-$ to the sequence \eqref{Eq-represent},
we obtain an exact sequence of left $A$-modules
\begin{eqnarray}\label{Eq-Asequence}
\Sigma^{\dag}\tensor{R}F_i \longrightarrow \Sigma^{\dag}\tensor{R}\Sigma_i \longrightarrow 0.
\end{eqnarray}
Since there is an $A$-module  isomorphism $$
\Sigma^{\dag}\tensor{R}F_i  \,\, \cong \,\, \oplus_{P\,\in\,
\cat{P}_i} P^*$$ for each $i \in I$, we deduce from the sequence
\eqref{Eq-Asequence} that $\Sigma^{\dag}\tensor{R}\Sigma_i$ is a
finitely generated $A$-module for each $i \in I$.

Since $\coring{C}$ is $R-\Sigma$-Galois, we have $\Sf{can}
:\Sigma^{\dag}\tensor{R}\Sigma \to \coring{C}$ is an isomorphism
of an $A$-corings. Therefore, the decomposition in
\eqref{Eq-comatrix} can be transferred via $\Sf{can}$ to a
decomposition of ${}_{\coring{C}}\coring{C}\,=\, \oplus_{i \,\in
\, I} E_i$ as a direct sum of subcomodules which are finitely
generated as left $A$-modules. Each one of the $E_i$'s is of
course a projective $A$-module as ${}_A\coring{C}$ is projective,
and this finishes the proof.
\end{proof}

According to \cite[Theorem 10]{Lin:1977} a coalgebra $C$ over a
field admits a decomposition as a direct sum of finite-dimensional
left subcomodules if and only if $\rcomod{C}$ has a generating set
of finite-dimensional projective right comodules ($C$ is already
right semiperfect). Thus, in the coalgebra case, the condition on
$\Sigma$ in statement $(ii)$ of Theorem \ref{TeoremaII}, namely
that ${}_R\Sigma$ is a direct sum of finitely generated
$R$--submodules, may be deleted. We do not know if this is also
the case for corings over a general ring (for Quasi-Frobenius
ground rings, the answer is positive \cite[Theorem
3.5]{ElKaoutit/Gomez:2005}). Note that the additional condition on
the left $R$--module structure of $\Sigma$ cannot be avoided in
statement $(iii)$ even in the coalgebra case, as there is a
generating set for $\rcomod{C}$ of finite-dimensional comodules
for any coalgebra $C$ over a field. This is also the case for
statement $(iv)$, since any coalgebra is $R-\Sigma$--Galois (by
\cite[Theorem 4.8]{ElKaoutit/Gomez:2004b} and \cite[diagram
(5.12)]{ElKaoutit/Gomez:2004b}).

As an application of Theorem \ref{TeoremaII}, we obtain the
following remarkable fact concerning to the existence of enough
projectives for categories of entwined modules. For the definition
of an entwining structure and a discussion of their properties and
their relationships with corings, we refer to
\cite{Brzezinski/Wisbauer:2003}.

\begin{corollary}\label{entgenerators}
Let $\Psi : A \tensor{} C \rightarrow C \tensor{} A$ be an
entwining structure between an algebra $A$ and a coalgebra $C$
over a commutative ring $K$. If $C$ admits a decomposition as a
direct sum of left subcomodules profinite as $K$--modules (e.g.,
if $K$ is a field and $C$ is right semiperfect), then the category
of right entwined modules $\entmod{A}{C}$ has a generating set of
small projective objects.
\end{corollary}
\begin{proof}
By \cite[Proposition 2.2]{Brzezinski:2002}, $A \otimes C$ is
endowed with the structure of an $A$--coring such that
$\rcomod{\coring{C}}$ is isomorphic to the category of right
entwined $A-C$--modules. The comultiplication on $A \otimes C$ is
given by
\[
\xymatrix{A \otimes C \ar^-{A \otimes \Delta_C}[rr] & & A \otimes C \otimes C \cong A \otimes C \tensor{A} A \otimes C}
\]
so that every decomposition of $C$ as direct sum of left
$C$--comodules leads to such a decomposition of the $A$--coring $A
\otimes C$. Obviously, if the direct summands in $C$ are profinite
$K$--modules, then the corresponding direct summands of $A
\tensor{} C$ are profinite as left modules over $A$. Now, the
Corollary is a consequence of Theorem \ref{TeoremaII}.
\end{proof}

Right semiperfect coalgebras over a field are characterized by the
fact that the rational functor $\Rat : \lmod{{}^*C} \rightarrow
\lmod{{}^*C}$ is exact \cite[Theorem 3.3]{Gomez/Nastasescu:1995}.
The exactness of the rational functor associated to a general
rational pairing in the context of corings has been considered
recently in \cite{ElKaoutit/Gomez:2008}. The rational functor
canonically associated to a coring becomes exact for the corings
characterized in Theorem \ref{TeoremaII}, as the following
Proposition shows. The exactness of the rational functor was given
in \cite{Brzezinski/Wisbauer:2003} (see Remark \ref{EKGTBW}
below).

\begin{proposition}\label{descexact}
Let $\coring{C}$ be an $A$--coring admitting a direct sum
decomposition $\coring{C} = \bigoplus_{E \in \cat{E}} E$, for a
family of subcomodules $\cat{E}$ of ${}_{\coring{C}}\coring{C}$
such that ${}_AE$ is profinite for every $E \in \cat{E}$. Consider
the canonical right rational pairing
$\TT=(\coring{C},\lend{\coring{C}}{\coring{C}}, \esc{-}{-})$, and
denote by $\fk{a} := \Rat^{\TT}\left(
\lend{\coring{C}}{\coring{C}}_{\lend{\coring{C}}{\coring{C}}}
\right)$ the rational ideal. Then
\begin{enumerate}[(a)]
\item $R$ is a right rational
$\lend{\coring{C}}{\coring{C}}$--module injected in  $\fk{a}$.
\item $\fk{a}$ is generated as a bimodule by $\{e_{E}\}_{E
\,\in\,\cat{E}}$, that is $\fk{a}= \sum_{E\,\in\,
\cat{E}}\lend{\coring{C}}{\coring{C}}
e_{E}\lend{\coring{C}}{\coring{C}}$, and $\fk{a}$ is a pure left
submodule of $\lend{\coring{C}}{\coring{C}}$.%
 \item The functor
$\Rat^{\TT}: \rmod{\lend{\coring{C}}{\coring{C}}} \rightarrow
\rmod{\lend{\coring{C}}{\coring{C}}}$ is exact.
\end{enumerate}
\end{proposition}
\begin{proof}
$(a)$ Follows directly from the isomorphisms \eqref{Eq-Iso} and
\eqref{Eq-EFR}. \\ %
\noindent $(b)$ From the isomorphisms \eqref{Eq-Iso} we get that
$\sum_{E\, \in\, \cat{E}} e_{E}\lend{\coring{C}}{\coring{C}}$ is a
rational right ideal, which implies that $\sum_{E\, \in\,
\cat{E}}\lend{\coring{C}}{\coring{C}}
e_{E}\lend{\coring{C}}{\coring{C}} \subseteq \fk{a}$. Thus, we
only need to check the reciprocal inclusion. So let $b \in
\fk{a}=\Rat^{\TT}\left(
\lend{\coring{C}}{\coring{C}}_{\lend{\coring{C}}{\coring{C}}}
\right)$ an arbitrary element with a right rational system of
parameters $\{(b_i,c_i)\}_{i\,=1,\cdots,n} \subseteq
\lend{\coring{C}}{\coring{C}} \times \coring{C}$. Let $\cat{E}'
\subset \cat{E}$ be a finite subset such that $c_i \in \oplus_{E
\in \cat{E}'}E$ for every $i = 1, \dots, n$, and take $e=\sum_{E\,
\in \cat{E}'}e_{E}$. One easily checks that $c_i\, e = c_i$, for
all $i=1,\cdots,n$. Therefore,
\begin{equation}\label{puro}
b\, e \,=\, \sum_{i} b_i \esc{c_i}{e}
\,=\, \sum_{i} b_i \varepsilon_{\coring{C}}(c_i \, e) \,=\, \sum_i
b_i \varepsilon(c_i) \,=\, b,
\end{equation}
which gives the needed inclusion.  Equation \eqref{puro} also
implies that $\fk{a}$ is a pure left
$\lend{\coring{C}}{\coring{C}}$--submodule of
$\lend{\coring{C}}{\coring{C}}$. \\ %
\noindent $(c)$ In view of $(b)$ and the equality
$\coring{C}\fk{a} = \coring{C}$, we can apply \cite[Theorem
1.2]{ElKaoutit/Gomez:2008} to get the exactness of the rational
functor.
\end{proof}

\begin{example}[{compare with \cite[Example 5.2]{Caenepeel/Iovanov:2005}}]
Let $\coring{C}$ be a cosemisimple $A$--coring.  By \cite[Theorem
3.1]{ElKaoutit/Gomez/Lobillo:2004c}, ${}_A\coring{C}$ and
$\coring{C}_A$ are projective modules. So we can consider its
right canonical rational pairing
$\TT=(\coring{C},\lend{\coring{C}}{\coring{C}},\esc{-}{-})$. The
structure Theorem of cosemisimple corings \cite[Theorem
4.4]{ElKaoutit/Gomez:2003} implies that $\coring{C}$ is a direct
sum of left $\coring{C}$-comodules where each of them is finitely
generated and projective left $A$-module. Thus for a cosemisimple
coring $\Rat^{\TT}$ has to be exact.
\end{example}

\begin{remark}\label{EKGTBW}
Under the anti-isomorphism of rings $\lend{\coring{C}}{\coring{C}}
\cong {}^*\coring{C}$ the rational pairing considered in
Proposition \ref{descexact} goes to a rational pairing
$(({}^*\coring{C})^{op}, \coring{C})$ that gives rise to the
``more usual'' rational functor $\Rat : \lmod{{}^*\coring{C}}
\rightarrow \lmod{{}^*\coring{C}}$ considered for instance in
\cite[Section 20]{Brzezinski/Wisbauer:2003}, where its exactness
was studied. In particular, one deduces from \cite[20.8,
20.12]{Brzezinski/Wisbauer:2003} that $\Rat :
\lmod{{}^*\coring{C}} \rightarrow \lmod{{}^*\coring{C}}$ is exact
whenever ${}_\coring{C}\coring{C}$ admits a direct decomposition
as assumed in Proposition \ref{descexact}. Their arguments run on
a different road than ours.
\end{remark}

\section{Local corings and Semiperfect corings}\label{sperfect}

Let $\coring{C}$ be an $A$--coring such that its category
$\rcomod{\coring{C}}$ of all right $\coring{C}$--comodules is a
Grothen\-dieck category. The coring $\coring{C}$ is said to be
\emph{right semiperfect coring} if every finitely generated right
$\coring{C}$--comodule has a projective cover. That is,
$\coring{C}$ is right semiperfect if and only if (by definition)
$\rcomod{\coring{C}}$ is a Grothendieck semiperfect category in
the sense of M. Harada \cite[Section 3, page 334]{Harada:1973a}.
Following  \cite[page 347]{Mares:1963} and \cite[Section 1, page
330]{Harada:1973a}, a right $\coring{C}$--comodule $P$ is said to
be a \emph{semiperfect right comodule} if $P_{\coring{C}}$ is a
projective comodule and every factor comodule of $P$ has a
projective cover. A right $\coring{C}$--comodule $P$ is said to be
a \emph{completely indecomposable comodule} if its endomorphisms
ring $\rend{\coring{C}}{P}$ is a local ring (i.e., its quotient by
the Jacobson radical is a division ring).

Recall from \cite[Corollary 1]{Harada:1973a} that if $\cat{A}$ is
a locally finitely generated Grothendieck category, then $\cat{A}$
is semiperfect if and only if $\cat{A}$ has a generating set of
completely indecomposable projective objects. Of course, this
result can be applied in particular to $\rcomod{\coring{C}}$,
whenever it is a locally finitely generated Grothendieck category.

\begin{remark}
Assume that $\rcomod{\coring{C}}$ is a Grothendieck category. It
seems to be an open question if it is locally finitely generated.
In the case when ${}_A\coring{C}$ is locally projective in the
sense of \cite{Zimmermann-Huisgen:1976}, then
$\rcomod{\coring{C}}$ is isomorphic to the category of all
rational left ${}^*\coring{C}$--modules (by \cite[Lemma
1.29]{Abuhlail:2003a} and \cite[Theorem
2.6']{ElKaoutit/Gomez/Lobillo:2004c}) and, therefore, the set of
all cyclic rational left ${}^*\coring{C}$--comodules generates the
category of right $\coring{C}$--comodules. In fact, in this case,
a right $\coring{C}$--comodule is finitely generated if and only
if it is a finitely generated as a right $A$--module (see, e.g.,
\cite[Lemma 2.2]{ElKaoutit/Gomez:2005}). Thus,
$\rcomod{\coring{C}}$ has a generating set of comodules that are
finitely generated as right $A$--modules, whenever
${}_A\coring{C}$ is locally projective.
\end{remark}

The following theorem is a consequence of \cite[Corollary
2]{Harada:1973a}, \cite[Theorem 3]{Harada:1973a} and \cite[Theorem
5.7]{ElKaoutit/Gomez:2004b} (see also \cite[Theorem
6.2]{Gomez/Vercruysse:2007}).

\begin{theorem}\label{semiperfectos}
Let $\coring{C}$ be an $A$--coring and $\cat{P}$ a set of right
$\coring{C}$--comodules. The following statements are equivalent.
\begin{enumerate}[(i)]
\item ${}_A\coring{C}$ is a flat module and $\cat{P}$ is a
    generating set of small completely indecomposable
    projective comodules for $\rcomod{\coring{C}}$; \item
    ${}_A\coring{C}$ is a flat module, $\rcomod{\coring{C}}$
    has a generating set consisting of finitely generated
    objects, $\coring{C}$ is a right semiperfect $A$--coring
    and $\cat{P}$ contains a set of representatives of all
    semiperfect completely indecomposable right
    $\coring{C}$--comodules;
    \item each comodule in $\cat{P}$ is finitely generated and
        projective as a right $A$--module, $\coring{C}$ is
        $(R,\Sigma)$--Galois, where $\Sigma = \bigoplus_{P \in
        \cat{P}}P$ and $R = \bigoplus_{P,Q \in
        \cat{P}}\hom{\coring{C}}{P}{Q}$, ${}_R\Sigma$ is
        faithfully flat, and $\rend{\coring{C}}{P}$ is a local
        ring for every $P \in \cat{P}$.
\end{enumerate}
\end{theorem}
\begin{proof}
$(i) \Rightarrow (iii)$. This is a consequence of \cite[Theorem
5.7]{ElKaoutit/Gomez:2004b}.\\
 $(iii) \Rightarrow (ii)$. By
\cite[Theorem 5.7]{ElKaoutit/Gomez:2004b}, ${}_A\coring{C}$ is a
flat module and $-\tensor{R}\Sigma: \rmod{R} \rightarrow
\rcomod{\coring{C}}$ is an equivalence of categories. From
\eqref{sufidem} we get the decomposition $R = \bigoplus_{P \in
\cat{P}} u_PR$, and, since $\rend{R}{u_PR} \cong u_PRu_P \cong
\rend{\coring{C}}{P}$ for every $P \in \cat{P}$, we deduce that
$\{ u_PR : P \in \cat{P} \}$ is a generating set of completely
indecomposable finitely generated projective objects for
$\rmod{R}$. Therefore, $\{ u_PR \tensor{R} \Sigma : P \in \cat{P}
\}$ becomes a generating set of completely indecomposable finitely
generated projective objects for $\rcomod{\coring{C}}$. Now, $u_PR
\tensor{R} \Sigma \cong P$ as a right
$\rcomod{\coring{C}}$--comodule, for every $P \in \cat{P}$. This
implies that $\cat{P}$ is a generating set of finitely generated
completely indecomposable projective comodules. By \cite[Corollary
1]{Harada:1973a}, $\coring{C}$ is right semiperfect. Now, if $U$
is any completely indecomposable projective right
$\rcomod{\coring{C}}$, then $U$ is an epimorphic image of a direct
sum of comodules in $\cat{P}$. We get thus that $U$ is isomorphic
to a direct summand of a sum of comodules in $\cat{P}$. By the
Krull-Schmidt-Azumaya theorem, we obtain that $U$ is isomorphic to
one comodule in $\cat{P}$.
 \\
$(ii) \Rightarrow (i)$ From the proof of \cite[Theorem
3]{Harada:1973a}, it follows that there is a set of finitely
generated projective completely indecomposable generators of
$\rcomod{\coring{C}}$. This obviously implies that $\cat{P}$ is a
generating set of projective finitely generated indecomposable
comodules.
\end{proof}

Let us record some information deduced from Theorem
\ref{TeoremaII} as a kind of ``Structure Theorem'' of right
semiperfect $A$--corings.

\begin{corollary}\label{semiperfectcorings}
Assume ${}_A\coring{C}$ to be locally projective. Then
$\coring{C}$ is right semiperfect if and only if there is a set
$\cat{P}$ of right $\coring{C}$--comodules such that
\begin{enumerate}
\item Every $P \in \cat{P}$ is profinite as a right
    $A$--module,
\item $\coring{C}$ is $R-\Sigma$--Galois, where $\Sigma =
    \bigoplus_{P \in \cat{P}} P$ and $R = \bigoplus_{P,Q \in
    \cat{P}}\hom{\coring{C}}{P}{Q}$,
\item ${}_R\Sigma$ is faithfully flat, and
\item $\rend{\coring{C}}{P}$ is a local ring, for every $P \in
    \cat{P}$.
\end{enumerate}
\end{corollary}

As mentioned before, a coalgebra over a field is right semiperfect
if and only if its rational functor is exact. This
characterization has been extended to the case of corings over
QuasiFrobenius rings, see \cite[Theorem
4.2]{ElKaoutit/Gomez:2005}, \cite[Theorem
4.3]{Caenepeel/Iovanov:2005},
\cite[20.11]{Brzezinski/Wisbauer:2003}. For a general ground ring
$A$, we obtain the following result.

\begin{corollary}\label{semiperfet-rat}
Let $\coring{C}$ be a right semiperfect $A$--coring such that
${}_A\coring{C}$ is a locally projective module. Consider the
canonical right rational pairing
$\T=(\coring{C},\lend{\coring{C}}{\coring{C}},\esc{-}{-})$ and let
$\Sigma$ and $R$ be as in Corollary \ref{semiperfectcorings}.
Assume that ${}_R\Sigma$ is projective and its Jacobson radical
$\mathrm{Rad}\left({}_R\Sigma\right)$ is a superfluous
$R$--submodule of ${}_R\Sigma$. Then the rational functor
$\Rat^{\T}: \rmod{\lend{\coring{C}}{\coring{C}}} \rightarrow
\rmod{\lend{\coring{C}}{\coring{C}}}$ is exact.
\end{corollary}
\begin{proof}
By Corollary \ref{semiperfectcorings}, $R = \bigoplus_{P \in
P}u_PR$, where $\{ u_P : P \in \cat{P} \}$ is a set of ortogonal
local idempotents. By \cite[Theorem 2]{Harada:1973b}, $\{ Ru_P : P
\in \cat{P} \}$ becomes a generating set of completely
indecomposable (finitely generated) projective objects for
$\lmod{R}$. Since ${}_R\Sigma$ is projective, it must be a direct
summand of a direct sum of left $R$--modules of the for $Ru_P$.
Now, the non-unital version of \cite[Theorem 5.5]{Mares:1963}
gives that ${}_R\Sigma$ has to be a direct sum modules of the form
$Ru_P$. By Theorem \ref{TeoremaII}, ${}_{\coring{C}}\coring{C}$ is
a direct sum of comodules, profinite as left $A$--modules.
Proposition \ref{descexact} gives the exactness of
$\Rat^{\cat{T}}$.
\end{proof}

Next, we deal with the case where $\cat{P}$ contains only one
comodule. This naturally leads to the notion of a local coring.

\begin{definition}
A right semiperfect $A$--coring $\coring{C}$ whose category of
right comodules has a unique (up to isomorphism) type of
completely indecomposable semiperfect comodule, will be called a
\emph{right local coring}.
\end{definition}

The following is the version of Theorem \ref{semiperfectos} for
the case where $\cat{P}$ is a singleton.

\begin{corollary}\label{locales}
Let $\coring{C}$ be an $A$--coring and
$(\Sigma_{\coring{C}},\rho_{\Sigma})$ a right
$\coring{C}$--comodule with endomorphism ring
$T=\rend{\coring{C}}{\Sigma}$. The following statements are
equivalent
\begin{enumerate}[(i)]
\item ${}_A\coring{C}$ is a flat module and
    $\Sigma_{\coring{C}}$ is a small completely indecomposable
    projective generator in the category of right comodules
    $\rcomod{\coring{C}}$;%
    \item ${}_A\coring{C}$ is a flat module, the category
        $\rcomod{\coring{C}}$ is locally finitely generated,
        $\coring{C}$ is a right local $A$--coring and
        $\Sigma_{\coring{C}}$ the unique semiperfect
        completely indecomposable right
    comodule (up to isomorphism);%
    \item $\Sigma_A$ is a finitely generated and projective
        module, $\coring{C}$ is $T-\Sigma$--Galois,
        ${}_T\Sigma$ is a faithfully flat module, and $T$ is a
        local ring.
\end{enumerate}
\end{corollary}

Our last proposition gives a connection between right local
corings and a simple co\-semi\-simple corings. Every cosemisimple
coring is a right and left semiperfect coring by \cite[Theorem
3.1]{ElKaoutit/Gomez/Lobillo:2004c}. A \emph{simple cosemisimple
coring} is a cosemisimple coring with one type of simple right
comodule, or equivalently, left comodule. The structure of these
corings was given in \cite[Theorem 4.3]{ElKaoutit/Gomez:2003}.

\begin{proposition}\label{simples-locales}
Let $\Sigma_A$ be a non zero finitely generated and projective
right $A$--module and $T \subseteq S=\rend{A}{\Sigma}$ a local
subring of its endomorphism ring such that ${}_T\Sigma$ is a
faithfully flat module. Consider the comatrix $A$-coring
$\coring{C}:=\rcomatrix{T}{\Sigma}$, and denote by $\mathrm{J}$
the Jacobson radical of $T$, and by $D:=T/\mathrm{J}$ the division
factor ring. Then we have
\begin{enumerate}[(i)]
\item $\Sigma/\mathrm{J}\Sigma$ admits the structure of a
    simple right $(\rcomatrix{T}{\Sigma})$--comodule whose
    endomorphism ring is isomorphic to the division ring $D$,
    that is $\rend{\coring{C}}{(\Sigma/\mathrm{J}\Sigma)}
    \cong D$. \item If $(\Sigma/\mathrm{J} \Sigma)_A$ is a
    projective module, then the canonical map
$$ \xymatrix@R=0pt{ \can{\Sigma/\mathrm{J}\Sigma}:
\left( \Sigma/ \mathrm{J}\Sigma\right)^* \tensor{D} \left(
\Sigma/ \mathrm{J}\Sigma\right) \ar@{->}[r] &
\rcomatrix{T}{\Sigma} }$$ is a monomorphism of $A$-corings
with domain a simple cosemisimple coring. \item There is a
short exact sequence of $\rcomatrix{T}{\Sigma}$--bicomodules
$$\xymatrix{0 \ar@{->}[r] & \Sigma^* \tensor{T} (\mathrm{J}\Sigma)
\ar@{->}[r] & \rcomatrix{T}{\Sigma} \ar@{->}[r] & \Sigma^*
\tensor{T} D\tensor{T} \Sigma \ar@{->}[r] & 0 }$$ whose
cokernel is an $A$--coring without counit.\item ${}_TS$ is a
faithfully flat module and the canonical Sweedler's
$S$--coring $S\tensor{T}S$ is a right local $S$--coring with
$S_{S\tensor{T}S}$ the unique (up to isomorphism) semiperfect
completely indecomposable comodule.
\end{enumerate}
\end{proposition}
\begin{proof}
We first make some useful observations. Recall that $\Sigma$ is a
right $\coring{C}$--comodule with coaction $\rho_{\Sigma}(u) =
\sum_i u_i \tensor{A} u_i^* \tensor{T} u$, for every $u \in
\Sigma$, where $\{(u_i,u_i^*)\}_i \subseteq \Sigma \times
\Sigma^*$ is any finite dual basis for $\Sigma_A$, while
$\Sigma^*$ is a left $\coring{C}$-comodule with coaction
$\lambda_{\Sigma^*}(u^*) = \sum_i u^* \tensor{T} u_i \tensor{A}
u_i^*$, for every $u^* \in \Sigma^*$. It is clear that
$\rho_{\Sigma}$ is left $T$-linear while $\lambda_{\Sigma^*}$ is
right $T$-linear. Since ${}_T\Sigma$ is a faithfully flat module,
\cite[Theorem 3.10]{ElKaoutit/Gomez:2003} implies that $T=
\rend{\coring{C}}{\Sigma}$ and $-\tensor{T}\Sigma: \rmod{T}
\rightarrow \rcomod{\coring{C}}$ is an equivalence of categories
with inverse the cotensor functor $-\cotensor{\coring{C}}
\Sigma^*: \rcomod{\coring{C}} \rightarrow \rmod{T}$. This means,
by Corollary \ref{locales}, that $\coring{C}$ is a right local
$A$-coring with $\Sigma_{\coring{C}}$ the unique type of
semiperfect completely indecomposable right comodule (up to
isomorphisms).

\noindent $(i)$ The right $\coring{C}$-coaction of the right
$A$-module $\Sigma /\mathrm{J}\Sigma$ is given by the functor
$-\tensor{T}\Sigma$, since $\Sigma /\mathrm{J}\Sigma \cong D
\tensor{T} \Sigma$ is an isomorphism of right $A$-modules. This
comodule is the image under the equivalence $-\tensor{T}\Sigma$ of
the simple right $T$-module $D_{T}$, so it is necessarily a simple
right $\coring{C}$-comodule. Its endomorphism ring is computed
using the following isomorphisms:
\begin{multline*}
\mathrm{Hom}_{\coring{C}}\left(D \tensor{T}\Sigma \, ,\,
D\tensor{T}\Sigma \right) \,\, \cong \,\, {\rm Hom}_T
\lr{D,\,\hom{\coring{C}}{\Sigma}{D\tensor{T}\Sigma}}
\\ \,\, \cong \,\,
{\rm
Hom}_T\lr{D,\,\hom{T}{(\Sigma\cotensor{\coring{C}}\Sigma^*)}{D}}
\,\, \cong \,\, \hom{T}{D}{D},
\end{multline*}
where we have used the adjunction $-\cotensor{\coring{C}}\Sigma
\dashv -\tensor{T}\Sigma$ in the second isomorphism and the right
$T$-linear isomorphism $\Sigma \cotensor{\coring{C}}\Sigma^* \cong
T$ in the last one.

\noindent $(ii)$ Since $(\Sigma/\mathrm{J}\Sigma)_A$ is a finitely
generated and projective module, we can construct the canonical
map $\can{\Sigma/\mathrm{J}\Sigma}$. This map is a monomorphism by
$(i)$ and \cite[Theorem 3.1]{Brzezinski:2005}. Finally,
\cite[Proposition 4.2]{ElKaoutit/Gomez:2003} implies that $\left(
\Sigma/ \mathrm{J}\Sigma\right)^* \tensor{D} \left( \Sigma/
\mathrm{J}\Sigma\right)$ is simple cosemisimple $A$-coring, since
$D$ is already is division ring.

\noindent $(iii)$ The stated sequence is obtained by applying
$-\tensor{T}\Sigma$ to the following short exact sequence of right
$T$-modules
$$\xymatrix{0 \ar@{->}[r] & \Sigma^*\mathrm{J} \ar@{->}[r] &
\Sigma^* \ar@{->}[r] & \Sigma^* \tensor{T} D \ar@{->}[r] & 0. }$$
This gives in fact an exact sequence of right
$\coring{C}$-comodules; the $\coring{C}$-bicomodule structure is
then completed by the compatible left $\coring{C}$-coaction of
${}_{\coring{C}}\Sigma^*$. The $\coring{C}$-bicomodule
$\Sigma^*\tensor{T} D \tensor{T} \Sigma$ admits a coassociative
comultiplication given by
$$\xymatrix@R=0pt{ \Delta: \Sigma^* \tensor{T} D \tensor{T} \Sigma
\ar@{->}[r] & \left( \Sigma^* \tensor{T} D \tensor{T} \Sigma
\right)
\tensor{A} \left(\Sigma^* \tensor{T} D\tensor{T} \Sigma \right) \\
u^* \tensor{T}d \tensor{T} v \ar@{|->}[r] & \sum_i u^*\tensor{T} d
\tensor{T} u_i \tensor{A} u_i^* \tensor{T}1 \tensor{T} v. }$$

\noindent $(iv)$ This is an immediate consequence of \cite[Theorem
3.10]{ElKaoutit/Gomez:2003} and Corollary \ref{locales}.
\end{proof}

\providecommand{\bysame}{\leavevmode\hbox
to3em{\hrulefill}\thinspace}
\providecommand{\MR}{\relax\ifhmode\unskip\space\fi MR }
\providecommand{\MRhref}[2]{%
  \href{http://www.ams.org/mathscinet-getitem?mr=#1}{#2}
} \providecommand{\href}[2]{#2}

\end{document}